\documentclass{article}

\usepackage{arxiv}

\usepackage{hyperref}       % hyperlinks
\usepackage{url}            % simple URL typesetting
\usepackage{booktabs}       % professional-quality tables
\usepackage{amsmath}        % blackboard math symbols
\usepackage{amsfonts}       % blackboard math symbols
\usepackage{nicefrac}       % compact symbols for 1/2, etc.
\usepackage{microtype}      % microtypography
\usepackage{cases}          % number cases
\usepackage{graphicx}
\graphicspath{ {./images/} }

\usepackage{bm}
\usepackage{setspace}
\usepackage{tabu}

\usepackage[sorting=none]{biblatex}
\bibliography{references}

\usepackage{caption}
\DeclareCaptionLabelFormat{Appendix}{A.#2}

\def\floor#1{\lfloor #1 \rfloor}
\def\ceil#1{\lceil #1 \rceil}

\newcommand{\divides}{\mid}

\hypersetup{
    colorlinks=true,
}

\title{Combined Sieve Algorithm for Prime Gaps}

\author{
 Seth Troisi \\
 \texttt{sethtroisi@google.com}
}
\begin{document}

\maketitle
\begin{abstract}
A new Combined Sieve algorithm is presented with cost proportional to the number of enumerated factors over a series of intervals. This algorithm achieves a significant speedup, over a traditional sieve, when handling many ($[10^4,10^7]$) intervals concurrently. The speedup comes from a space-time tradeoff and a novel solution to a modular equation. In real world tests, this new algorithm regularly runs 10,000x faster. This faster sieve paired with higher sieving limits eliminates more composites and accelerates the search for large prime gaps by 30-70\%. During the development and testing of this new algorithm, two top-10 record merit prime gaps were discovered.
\end{abstract}

\keywords{primes, prime differences, prime gaps, parallel sieve}

\section{Introduction}

\iffalse

 Todo:
\begin{itemize}
    * on Xeon K = 807#/7#, M = 200.5m => 45,800,000 in 8679m => 87.95/s
    * Graph of factors / second
\end{itemize}

\fi

A prime gap is the difference between successive primes. A record is kept of the smallest known prime for each gap. For example the smallest prime followed by a gap of 4 is 7; a gap of 6 first follows 23; and a gap of 96 first follows 360,653. Until his death in 2019, Thomas R. Nicely maintained a list of record prime gaps, first in print and later on his faculty website \cite{Nicely1999NewMP}. Since 2019, a group from MersenneForum.org, including the author, has maintained the records on GitHub\cite{primegaplist}.

The search for prime gaps precedes computers but the digital era accelerated the search. Computer searches for prime gaps start circa 1930 when D.H. Lehmer constructed a physical computer from bicycle chains, The Lehmer Sieve, and used it to construct a table of the primes less than 37 million \cite{Selfridge1959TablesCT}. In 1961 Baker, Gruenberger, And Armerding added to the search in \citetitle{gruenberger1961statistics} \cite{gruenberger1961statistics} Both listed first occurrences of prime gaps. Lander and Parkin using a CDC 6600 extended the search and published a paper exclusively on prime differences in 1967 \cite{Lander1967OnFA}. At least a dozen papers from as many authors have followed since.

This paper outlines a new algorithm that speeds up searches for large gaps by combining the sieving of many nonadjacent but regularly spaced intervals at a substantial speedup. With a faster sieve more small factors can be found reducing the number of candidates tested per interval. Ready-to-run code, benchmarks, and validation can be found on GitHub\footnote{\url{https://github.com/sethtroisi/prime-gap}}.

\subsection{New Prime Gaps with Merit 38.048 \& 37.051}

Large gaps become more common as the density of primes decrease. The merit of a prime gap, $g_k$, following the prime, $p_k$, is defined as $\frac{g_k}{\log(p_k)}$. It is the ratio of the gap to the ``average" gap near that point. Appendix \ref{figure:recordmerit} shows the merit of each record prime for gaps under $100,000$. Gaps with merit $> 35$ are exceedingly rare, only 30 are known, so finding two using this new algorithm is both an accomplishment and involved quite a bit of luck. The two new prime gaps are:

A prime gap of $\textbf{35,308}$ follows $\bm{100054841 \frac{953\#}{210} - 9670}$ with merit = $\textbf{38.048}$ this is the \textbf{4th} largest known merit\footnote{\url{https://primegap-list-project.github.io/lists/top20-overall-merits/}}.

A prime gap of $\textbf{36,702}$ follows $\bm{164065661 \frac{1019\#}{210} - 18566}$ with merit = $\textbf{37.051}$ this is the \textbf{9th} largest know merit.

\section{Searching for Prime Gaps}
\label{sec:headings}

An exhaustive search for the first occurring prime gap has been performed up to $2^{64}$\cite{primegapfullyanalyzed}. Thomas Nicely explains the search process

\begin{quote}
     No general method more sophisticated than an exhaustive search is known for the determination of first occurrences and maximal prime gaps. As in the present study, this is most efficiently done by sieving successive blocks of positive integers for primes, recording the successive differences, and thus determining directly the first occurrences and maximal gaps. \cite{Nicely1999NewMP} 
\end{quote}

The largest known first occurring prime gap is the gap of $1550$ following $18361375334787046697$ with merit $\frac{1550}{\log(1.836e19)} = \textbf{34.944}$ found by Bertil Nyman in 2014 and independently confirmed in 2018\cite{primegaplist}. A history of the upper bound for exhaustive searches is found on Thomas Nicely's page \cite{NicelyFaculty} and an updated version on GitHub \cite{primegapfullyanalyzed}.

\paragraph{Larger Gaps}
The natural approach to search for larger gaps is similar. First, an interval around a test number, $N$, is sieved using small primes to mark numbers with small factors as composite. Then a primality test is run on the remaining numbers to find the ''surrounding" primes, the next and previous prime. The prime gap containing $N$ is the difference between the two surrounding primes. An existing respected program \texttt{PGsurround.pl}\cite{PGSurround} use this approach with the additional optimization of searching only $N$ of a certain forms known to produce larger gaps.

The speed a program using this approach can find the primes surrounding $N$ is determined by the speed of the primality test and how many composite numbers are removed by the sieve. Generally less than 2\% of execution time is spent on the sieve stage; the Combined Sieve doesn't change this dynamic, even with the deepest sieves, sieving up to $10^{11}$, and smallest $N$, $N \sim 3,000$ bits, sieving takes less than $10\%$ of amortized runtime (see Appendix \ref{appendix:runtime}).

\paragraph{Primorials}

It's possible to construct arbitrarily large gaps by looking near primorials.

Let $K = P\# = \displaystyle\prod_{p \in \mathbb{P}}^{P} p\ $

\begin{align*}
\begin{array}{rll}
    2 \divides & (K + 2) \implies & K + 2 \notin \mathbb{P} \\
    3 \divides & (K + 3) \implies & K + 3 \notin \mathbb{P} \\
    2 \divides & (K + 4) \implies & K + 4 \notin \mathbb{P} \\
    & ... & \\
    & (K + X) & \\
    & ... & \\
    P \divides & (K + P) \implies & K + P \notin \mathbb{P}
\end{array}
\end{align*}

Thus there are no primes in the interval $[K+2, K+P]$ and the gap from $prev\_prime(K+2)$ to $next\_prime(K+2)$ is $\ge P$. If both $K\pm 1$ are composite (by far the average case) then the gap is $\geq 2P+3$. Note that as the gap grows larger so does $\log(K)$ so merit doesn't grow arbitrarily large.

The average gap can be increased by choose a small high composite divisor $d$ (e.g. $11\#$) and looking near $K = \frac{P\#}{d}$. To illustrate why, consider $K = \frac{P\#}{2}$ all odd $X$, even those $> P$, are multiples of two ($K$ is odd so $K + X$ is even) and multiples of small primes eliminate all even $X \le 2P$ that aren't powers of two.

Because of the way primorials align their prime factors a greater share of nearby numbers are composite, meaning larger gaps. For this reason nearly all of the records ($\approx 97\%$ of the 95,000+ record gaps with gap $> 2000$) are of the form $m\frac{P\#}{d} - a$, $m$ being a small multiple, $a$ being the offset to the first prime below $m\frac{P\#}{d}$. 

\subsection{Sieving}

For each $N$, an interval (typically a few thousand numbers corresponding to average or optimistic prime gap) is sieved around $N$ and numbers without any factors, heading away from $N$, are tested for primality till the surrounding primes are found.

The traditional sieve calculates, for all small primes $p$, $N \equiv r \bmod p$ and uses that to mark $N - r \pm i \times  p$ as composite. E.g. $N \equiv 5 \mod 73 \implies 73 \divides N-5,\ \ 73 \divides N-5-73 \text{,\ \ and\ \ } 73 \divides N-5+73$.

\subsection{Probabilistic Primality Testing}

The best known deterministic primality tests are much slower than their probabilistic counterparts with respectively runtime of $\tilde{\mathcal{O}}(\log(n)^6)$ and $\tilde{\mathcal{O}}(\log(n)^2)$. As the probability of a composite being marked as prime can be reduced to an arbitrary level (BPSW has no know counterexamples) and no primes is ever marked as composite, it's advantageous to use a probabilistic primality tests. Notable gaps, with high merit, often have their endpoints later certified as prime with a deterministic test.

GMP version 6.2.0 \cite{gmplib} \href{https://gmplib.org/manual/Number-Theoretic-Functions#Number-Theoretic-Functions}{\texttt{mpz\_probab\_prime\_p}} uses a BPSW primality test. GMP is well optimized for numbers of all sizes and is a natural choice for primality testing\footnote{Some of the inspiration for this paper came from analysis done while improving GMP's \texttt{next\_prime} code}. Some testing was done with PrimeForm/GW which is faster for large numbers (1000+ digits) but more difficult to work with.

\section{Combined Sieve Algorithm}

For reasons described above, the majority of prime gap searches are centered around numbers of the form $N = m \frac{P\#}{d}$. Often millions of values of $m$ are tested for the same $K = \frac{P\#}{d}$. The Combined Sieve algorithm sieves $M$ sequential intervals simultaneously at a dramatically speedup to sieving $M$ intervals serially. The algorithm works on any arithmetic sequence but for the purpose of prime gap it's beneficial to focus on numbers of this deficient primorial form.

\begin{figure}[ht]
    \center
    \includegraphics[scale=0.62]{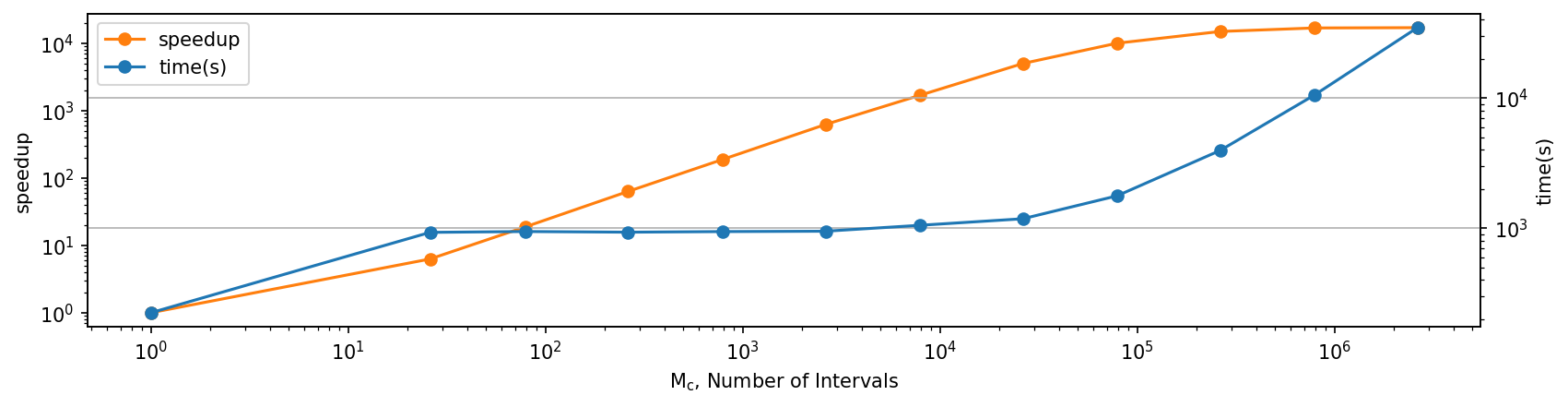}
    \caption{Speedup for various values of $M_c$ with $K= \frac{1511\#}{7\#}, X=18,000,$ and $P_{limit} = 10^{11}$}
    \label{figure:speedup}
\end{figure}

Note that $m$'s that share a factor with $d$ are skipped partially to avoid duplicate testing (e.g. $m=100, d=35$ would be a repeat of $m=20, d=7$) and partially because small $d$ tend to generate lower merit gaps. Figure \ref{figure:speedup} accounts for this by counting the number of $m$ coprime to $d$ less than $M$, $M_c$ = $| \{ m : 1 \le m \le M \land  (m, d) = 1 \} |$

\subsection{Space–Time Tradeoff}
\label{subsec:trick1}

By caching $K \bmod p$, $N \bmod p$ is reduced to a machine word sized modulo.

\begin{equation}
\begin{split}
    N \bmod p = m \frac{P\#}{d} \bmod p &= m K \bmod p \\
      &\equiv (\underbrace{(m \bmod p)}_\text{$<$64bits} \underbrace{(K \bmod p)}_\text{500-20,000 bits}) \bmod p
\end{split}
\end{equation}

After this trade-off only the first $m$ is slow for the other $M_c-1$ values of $m$ each $N \bmod p$ can be computed with a memory lookup, a 64 bit multiply and a 64 bit modulo irregardless of the size of $K$. This means the max speedup is $M_c$ since $K \bmod p$ is still calculated once for each $p$. This upper bound is easily seen in Figure \ref{figure:speedup} up to $M_c \sim 10^5$.

There are two easy ways to cache $K \bmod p$. $K \bmod p$ can be caching for all primes than two nested for-loops iterate over $m, p$ this requires storing the modulo for all primes in memory. The other method inverts the order and iterating over $p, m$ this requires storing the list of coprime $m$ and a very large boolean array in memory of the composite status of all offsets for all $M_c$ intervals.

\subsection{Large Prime Optimization}
\label{subsec:trick2}

This section details a complicated trick that allows the algorithm to avoid calculating $m K \bmod p$ for each $m$ and instead directly compute the set of $m$ where $p$ divides a number in the intervals around $N = m K$.

Assume $X$ numbers before and after $N$ are sieved for an interval of size $S=2X+1$, on average $\frac{S}{2}$ numbers will be divisible by $2$, $\frac{S}{3}$ by $3$, $\frac{S}{5}$ by $5$. In general, each prime $p$ will divide $\ceil{\frac{S}{p}}$ or $\floor{\frac{S}{p}}$ numbers in the interval around $N$. When $p \gg S$ most primes will not divide a number in the interval. Mertens' second theorem states

\begin{equation}\label{eq:meisselmertens}
    \sum_{p}^{n} \frac{1}{p} = \log(\log(n)) + M_m + o(1)    \\
\end{equation}
Where $M_m \approx \ 0.2615$ is the Meissel-Mertens constant.

\begin{spacing}{1.25}
With $X=100000, N=\frac{503\#}{210}$, all of the $17,984$ primes $2 \leq p \leq S$ divide one or more number in the interval. $\frac{24744}{60514} = 41\%$ of the  primes $S < p \leq 10^6$ divide a number in the interval. $\frac{81355}{50769036} =0.16\%$ of the primes $10^6 < p \leq 10^9$ divide a number in the interval. This means that $99.84\%$ of large primes for any given $m$ don't divide a number in the interval but still require processing. As the sieve limit increases, to $10^{10}$, $10^{11}$, $10^{12}$, ..., the fraction of primes dividing a number in the interval decreases further.
\end{spacing}

\subsubsection{Testing many intervals}

With a fixed prime $p$ and fixed $K$ and $1 \leq m \leq M$, sieving can be restated as trying to finding $m$'s where

\begin{equation}\label{eq:modeq}
\begin{split}
    -X \le\ &m K \le X \bmod p \\
    &\equiv \\
    0 \le \ &m K + X  \le 2X \bmod p \\
\end{split}
\end{equation}

For each $p$ larger than some threshold (naively $p > 10 S$ i.e. $p$ divides fewer than 1 in 10 intervals), $K \bmod p$ is calculated, and the smallest $m$ satisfying Equation \ref{eq:modeq} is found. $m  K \bmod p$ is computed and one of $\floor{\frac{m K}{p}}p$ or $\ceil{\frac{m K}{p}}p$ is marked as a composite. The algorithm continues finding larger $m$ and marking off composites till $m > M$.

Because we're only interested in $m$ which are coprime to $d$, work is wasted when $m$ shares a factor with $d$ and the search advances to the next $m$ without marking a composite. The smaller and more numerous the factors of $d$, the more unnecessary work. Very large $d$ are uncommon so a reasonable worst case is $d=17\#$ where $81.9\%$ of $m$ share a factor with $d$.

Each solution to Equation \ref{eq:modeq} (ignoring $m$ which share a factor with $d$ and the last solution with $m > M$) is another factor of a number in one of the interval. This is a very powerful transformation, by solving Equation \ref{eq:modeq} the Combined Sieve algorithm can enumerate $m$'s where $p$ divides a number in the interval $[m K - X, m K + X]$. Each prime divides an average $\frac{M S}{p}$ numbers because $p \gg S$ this is a large reduction from testing all $M_c$ intervals individually. With $S=60,001$ and $P \sim 10^{10}$ only 1 in $160,000$ $m$ will need to calculate $mK \bmod p$. 

Equation \ref{eq:meisselmertens} gives an estimate of the number of primes $p \leq P_{limit}$ that divide a number in each interval as $\approx S\ (\log(\log(P_{limit})) + o(1))$. This means that all prime factors $10 S \leq p \leq n$ for all $M$ intervals can be enumerated by solving Equation \ref{eq:modeq} $\approx M \times S (\log(\log(P_{limit})) - \log(\log(\text{10S}))) + \text{PrimePi}(P_{limit}) - \text{PrimePi}(10S)$ times\footnote{This is one of the many theoretical results verified during execution. The count of all $m$ tested over all $p$ generally matches to $\approx0.05\%$}.

\subsection{Solving Equation \ref{eq:modeq}}

\begin{equation}\label{eq:modfunc}
    \text{Let }G(p, K, L, R) = m\text{ be the smallest solution to }L \le m K \le R \bmod p
\end{equation}

\begin{subnumcases}{\label{eq:modsol} G(p,K,L,R) = }
        0                     & $K = 0\ \|\ L = 0$                        \label{eq:modsol_a}\\
        G(p, K \bmod p, L, R) & $K \ge p$                                 \label{eq:modsol_B}\\
        \ceil{L/K}            & $\ceil{L/K} \le \floor{R/K}$            \label{eq:modsol_C}\\
        G(p,p-K,p-R, p-L)     & $2K \ge p$                              \label{eq:modsol_D}\\
          \ceil{(L + \frac{p \times G(K, -p \bmod{K}, L \bmod{K}, R \bmod{K}))}{K}} &  \label{eq:modsol_E}
\end{subnumcases}

This solution assumes $p$ and $K$ are coprime which is guaranteed by the definition of $K$ as a primorial and $p$ a large prime (at minimum $p > S > L, R$). $p$ and $K$ coprime also guarantees the existence of a solution $(L \times K^{-1}) \bmod p$.

A handwavy proof that \ref{eq:modsol} returns the minimal answer can be constructed by inspecting \ref{eq:modfunc} with and without the $\bmod\ p$ condition. If Equation \ref{eq:modfunc} has a solution without considering $\bmod\ p$, the minimal solution is the first multiple of $K \geq L \implies m = \ceil{L/K}$ , otherwise $L \le K m - i p \le R$ for some $i \geq 1$.

\begin{align*}
\begin{array}{rlll}
    L           \le & \ m K - i  p      & \le R         & = \\
    L - m K     \le & -i \times p       & \le R - m K   & \\
                \\
        \multicolumn{4}{c}{\text{consider} \bmod{K}}    \\
                \\
    L           \le & i \times -p       & \le R \bmod{K}  \\
        \multicolumn{4}{c}{i = G(K, -p \bmod{K}, L \bmod{K}, R \bmod{K})} \\
                    & \implies                             \\
        \multicolumn{4}{c}{m = \ceil{\frac{L + p G(K, -p \bmod{K}, L \bmod{K}, R \bmod{K})}{K}}}
\end{array}
\end{align*}

\ref{eq:modsol_D} which reduces $K > \frac{p}{2}$ to $K < \frac{p}{2}$ is the result of some rearrangement. 

\begin{align*}
\begin{split}
    \text{Let }G(p, K, L, R) = m &\text{ be the smallest solution to } \\
    L \le  m K \le R \bmod{p}   & \equiv  \\
    -L \ge m (-K) \ge -R & \equiv \\
    -R \le m (-K) \le -L & \equiv \\
    \text{Which is equivalent to} & \\
    G(p, -K, -R, -L) & \equiv G(p, p-K, p-R, p-L)
\end{split}
\end{align*}

Equation \ref{eq:modsol} either directly returns $m$ (\ref{eq:modsol_a}, \ref{eq:modsol_C}) or reduces to a new sub-problem with $K < \frac{p}{2}$ \eqref{eq:modsol_D} or $p < K$ \eqref{eq:modsol_E}. It follows that $G$ is $\mathcal{O}(\log(p))$.

Given $p < 2^{64}$ and $K \bmod p$ is cached, this runs quickly. In practice it takes $2\log_2(\frac{p}{R-L}) \leq 2\log_2(\frac{10^{13}}{10000}) \leq 60 $ steps. On an i7-2600K with $p$ a 50 bit prime and $R-L=S=100,000$ an optimized implementation takes 100-200 ns\footnote{See \url{https://github.com/sethtroisi/prime-gap/blob/main/BENCHMARKS.md}}.

\section{Parameter Selection}

\subsection{Sieve Limit}

Mertens' third theorem estimates the fraction of numbers that don't have a factor less than $P_{limit}$\\
\begin{equation}\label{sieveeffic}
    \prod_{p}^{P_{limit}} \frac{p - 1}{p} \approx \frac{1}{\log(P_{limit}) e^\gamma}
\end{equation}

Where $\gamma \approx 0.577$ is the Euler–Mascheroni constant.

\begin{center}
\begin{tabu}{ | c c c c |}
    \hline
    $P_{limit}$ & \textbf{primes} &
    \textbf{theoretical} & \textbf{experimental} \\
    \hline
    $10^4$ & 1229     & 0.060960 & 0.06088      \rule{0pt}{1em}\\ 
    $10^5$ & 9592     & 0.048768 & 0.04875      \\ 
    $10^6$ & 78498    & 0.040640 & 0.040638     \\ 
    $10^8$ & 5761455  & 0.030480 & 0.030480     \\
    $4 \times 10^9$ & 189961812    & 0.025394 & 0.025394   \\
    $10^{10}$ & 455052511   & 0.024384 & 0.024384   \\
    \hline
    $10^{11}$  & $4.1 \times 10^9$    	& 0.022167 & \rule{0pt}{1.1em}\\
    $10^{12}$  & $3.8 \times 10^{10}$   & 0.020320 & \\
    $10^{13}$  & $3.5 \times 10^{11}$   & 0.018757 & \\
    \hline
\end{tabu}
\end{center}

Joined with the fact that a number near $N$ is prime $\approx\frac{1}{\log(n)}$ of the time, the expected number of primality tests that will need to be performed to find a prime after removing numbers with factors less than $P_{limit}$ is

\begin{equation}\label{eq:prptests}
    E(N, P_{limit}) \approx \frac{\log(N)}{\log(P_{limit}) e^\gamma}
\end{equation}

A sieve limit can be chosen such that the quasi-linear cost of calculating $K \bmod p$ and solving Equation \ref{eq:modeq} multiplied by the expected number of primality tests reduced is near the cost of running the avoided primality tests.

In practice \texttt{combined\_sieve} computes the derivative of Equation \ref{eq:prptests} at regular intervals giving an estimate of primality tests avoided per second of additional sieving. When that rate falls below the benchmarked rate of primality tests per second, \texttt{combined\_sieve} advises the user to stop sieving. This allows for better control of sieve limit based on CPU and memory characteristics.

A counter-intuitively observation is the optimal sieve limits can decrease as $N$ increases. The optimal limit strongly depends on $M$ and for large $N$ it doesn't make sense to sieve a lot of intervals when only a few can be tested per day.

Anecdotally setting $P_{limit}$ much higher than $5 S M$ is sub-optimal as the ratio of time spent performing $K \bmod p$ starts to dwarf the time spent solving Equation \ref{eq:modeq}. Some values used during exploratory searches are listed in Appendix \ref{appendix:sievelimits}.

\subsection{Sieve Interval}

For any $N$, if both of the surrounding primes is within $X$ of $N$ then time was wasted finding factors of numbers further away. If either of the surrounding primes is further away than $X$ additional traditional sieving will need to be done. Choosing $X$ in $[6 \log(N), 10 \log(N)]$ is reasonable. In practice \texttt{combined\_sieve} performs a small test sieve of a few $m$ and issues a warning if the estimated probability of the surrounding primes exceeding $X$ is greater than 0.5\%.

\section{Runtime}

Let

\begin{align*}
    S   =&\ 2X + 1 = \text{the sieve interval} \\
    M   =&\ \text{number of combined sieves to run simultaneously} \\
    M_c =&\ \text{number of $m$ coprime to $d$} \\
    C_{mult}   =&\ \text{small prime multiplier constant} \approx 10\frac{M}{M_c} \sim 50 \\
    P_{limit}  =&\ \text{sieve limit} \\
    |P| =&\ \text{PrimePi}(P_{limit}) \approx \frac{P_{limit}}{\log(P_{limit})}  \\
    |P_{small}| =&\ \text{PrimePi}(C_{mult} S) \approx \frac{C_{mult} S}{\log(C_{mult} S)} \\
    |P_{large}| =&\ |P| - |P_{small}|
\end{align*}

\textbf{Note:} This section ignores the small one-time cost of finding primes up to $P_{limit}$.

Sieving each coprime $m$ without the \textbf{Space-Time Tradeoff} and \textbf{Large Prime Optimization} involves computing $(m \times K) \bmod p$ then marking off an average of $\frac{S}{p}$ composites. $N \bmod p$ can be computed in $\log(N)$ time,

\begin{equation}
    |P| \log(N) + S (\log(\log(P_{limit})) + M_m)
\end{equation}

\begin{spacing}{1.1}
The Combined Sieve algorithm computes $(K \bmod p)$ once and reused this for each $m$. Small primes, which are likely to divide a number in the sieve interval (e.g. $\frac{p}{S} < C_{mult}$) are handled similarly to the traditional sieve but use cached $K \bmod p$ to calculating $N \bmod p$ with 64 bit integers. For larger primes the $(m \times K) \bmod p$ calculation is avoided for the vast majority of $m$ at the amortized cost of $\frac{M S}{p} + 1$ calls to Equation \ref{eq:modeq}. The cost for $M$ intervals is
\end{spacing}

\begin{multline}
    \underbrace{|P|\ \log(K)}_{K \bmod p \text{ for all primes}} +
    \underbrace{M_c S (\log\log(P_{limit}) + M_m))}_{\text{marking off all multiples}} + \\
    \underbrace{M_c |P_{small}|}_{m \times K \bmod p \text{ for small primes}} +
    \underbrace{(M S (\log\log(P_{limit}) - \log\log(C_{mult} S)) + |P_{large}|)}_{\text{calls to Equation \ref{eq:modeq}}} \times
    \underbrace{\log_2(\frac{P_{limit}}{S})}_{\text{cost to solve Equation \ref{eq:modeq}}}
\end{multline}

With

\begin{align*}
    K =& \frac{8887\#}{13\#} \\
    \log_{2}(K) =& 12,651 \ \ (\log(P\#) \approx \log(e^P) \approx P) \\[0.3em]
    M =& 10,000 \\
    M_c =& 1,917 \\
    P_{limit} =& 10^{10} \\
    |P| =& 455,052,511 \\[0.5em]
    & \text{ Optimizing the searching above and below for up to 10 merit} \\
    S =& 2X+1 = 2 \times 10 \times 9000 + 1 = 180,001 \\[0.5em]
    C_{eq\ref{eq:modeq}} =& 30\ \text{experimental constant for speed of Equation \ref{eq:modeq}} \\
    C_{mod} =& \frac{1}{16}\ \text{experimental constant for speed of modulo} \footnotemark
\end{align*}

\footnotetext{\label{fn:cycles}On a i7-2600K @3.4 GHz, GMP computes $N \bmod p$ in $\approx \frac{4}{64}\log_2(N)$ cycles, Equation \ref{eq:modeq} takes $\approx 30 \log_2(\frac{P}{X})$ cycles to solve.}

The experimental speedup of the Combined Sieve over the traditional sieve is

\begin{equation}
\begin{split}
    & |P| \log_2(K) C_{mod} + M_c S (\log(\log(P_{limit})) + M_m) + \\
    & M_c |P_{small}| + (M S (\log(\log(P_{limit})) - \log\log(C_{mult} S)) + |P_{large}|) C_{eq\ref{eq:modeq}} \log_2(\frac{P_{limit}}{S}) \\
    \approx\ & 455,052,511 \frac{12,651}{16} + 1,917 \times 180,001 (3.14 + 0.26) + \\
    & 1,917 \times 602,492 + (10,000 \times 180,001 \times (3.14 - 2.86) + 452,618,854) \times 30 \log_2(\frac{10^{10}}{400,001})   \\
    \approx\ & 3.60 \times 10^{11} + 1.17 \times 10^{9} + 1.15 \times 10^{9} + 5.24 \times 10^{11} \\
    \approx\ & 8.86 \times 10^{11} \\
    \textbf{vs} & \\
    & M_c (|P| \log_2(K) C_{mod} + S (\log(\log(P_{limit})) + M_m)) \\
    \approx\ & 1,917 \times (455,052,511 \frac{12,651}{16} + 180,001 (3.14 + 0.26)) \\
    \approx\ & 6.90 \times 10^{14} \\[0.3em]
    \frac{6.90 \times 10^{14}}{8.86 \times 10^{11}} =\ &\ 778 \text{x speedup}
\end{split}
\end{equation}

This result matches experimental observations of speedups. In general the speedup mostly depends on the magnitude of $M$ (The speedup is bounded by $M_c$, $1,917x$ in this example), moderately on $P_{limit}$, and a small amount on the magnitude of $N$.

\section{Results}

Traditionally sieving with primes twice as large, means finding half as many factors in the same amount of time. Using the Combined Sieve algorithm $K \bmod p$ takes the same amount of time irrespective of $p$ but larger primes require solving Equation \ref{eq:modeq} fewer times. As long as a reasonable percentage of time is spent solving Equation \ref{eq:modeq} vs calculating $K \bmod p$, the sieve speeds up as it processes larger primes. Anecdotally it takes less than twice as long to sieve to a 10x higher sieve limit (e.g. if it takes 2 hours to handle primes up to 50 million, all primes up to 500 million can be handle in less than 4 hours).

Use of the Combined Sieve allows for substantial higher sieve limits which reduces the number of primality tests that need to be perform. Appendix \ref{appendix:sievelimits} has estimates of new optimal sieve limits. Calculating the resulting speedups compared to the sieve limits used by \texttt{PGSurround.pl} the overall speedup from using the Combined Sieve is \textbf{30-70}\%.

\paragraph{Real world searching}\label{realworldsearching}

A number of complex relationships govern overall speed. To maximize the number of records found a wide search strategy is utilized, various combinations of $\frac{P\#}{d}$ are sieved with $M = 2,000$. For each $m_n$ an auxiliary program computes $Prob_n = P(\text{record gap} | m_n \times K, \text{sieved interval})$, for the top 5\% $\{Prob_n, m_n, K, \text{sieved interval}\}$ is stored in a database. Several background threads read from the database testing the highest likelihood intervals. If $\overline{\frac{Prob}{\text{time}}}$ for a $\frac{P\#}{d}$ is above average, \texttt{combined\_sieve} is run with the largest $M$ that fits in memory and a higher sieve limit. A number of additional metrics and stats are collected to verify assumptions and integrity during this process; see Appendix \ref{appendix:probrecord} for plots of theoretical and experimental distributions.

\subsection{Future work}

\texttt{combined\_sieve} contains a number of other significant speedups not discussed in this paper. Section \ref{realworldsearching} briefly describes an optimization involving directly computing the probability of a record gap for each interval using the remaining numbers in an interval. This can be utilized to only test the most likely 5-30\% of $m$ or to avoid testing the high sides if previous prime is close to $N$.

\paragraph{Acknowledgments}

The author would like to thank Craig Citro for advice on publication and Iulia Turc for regular check-ins on progress.

\printbibliography

\appendix
\section*{Appendix}\
\captionsetup{labelformat=Appendix}
\setcounter{table}{0}
\setcounter{figure}{2}

\begin{table}[ht]
    \caption{Experimental Sieve Limits}
    \label{appendix:sievelimits}
    \center
    \begin{tabu}{ | c c c c c c |}
        \hline
        $P$ & $\log(P\#)$ & \textbf{sieve limit} & \textbf{expected primality tests} & \texttt{PGSurround.pl} sieve limit & Theoretical Speedup \\
        \hline
         1511 & 1469	&  $250 \times 10^9$	& 31.5	& $5.0 \times 10^6$	& 70.2\%	\rule{0pt}{1em}\\
         2221 & 2054	&  $800 \times 10^9$	& 40.9	& $1.4 \times 10^7$	& 71.8\%	\\
         4441 & 4353	& $2000 \times 10^9$	& 86.3	& $1.3 \times 10^8$	& 51.4\%	\\
         5333 & 5228	& $2000 \times 10^9$	& 103.7	& $2.3 \times 10^8$	& 47.0\%	\\
         8887 & 8779	& $1000 \times 10^9$	& 178.4	& $1.1 \times 10^9$	& 32.5\%	\\
        \hline
    \end{tabu}
\end{table}

\begin{table}[ht]
    \caption{Combined Sieve Amortized Runtime (time per coprime m)}
    \label{appendix:runtime}
    \center
    \begin{tabu}{ | l l r r r r r |}
        \hline
        $M$ & $d$ & $X$ & $M$ & sieve limit & Combined Sieve(s) & Testing(s) \\
        \hline
        1511 & $7\# \times 1009$    & 16,000    & 10,000,000 &  $100 \times 10^9$   & 0.0112    & 0.134\rule{0pt}{0.85em}\\
        1511 & -                    & -         & -           &  $250 \times 10^9$   & 0.0129    & 0.127     \\
        2221 & $5\# \times 2027$    & 23,000    &  10,000,000 & $800 \times 10^9$   & 0.0252    & 0.568     \\
        4441 & $11\# \times 4001$   & 47,000    &  10,000,000 & $2000 \times 10^9$   & 0.0605    & 6.04      \\
        5333 & $11\# \times 5003$   & 65,000    &   1,000,000 & $2000 \times 10^9$   & 0.158     & 11.23     \\
        8887 & $17\# \times 8009$   & 100,000   &     100,000 & $1000 \times 10^9$   & 0.709     & 72.16     \\
        \hline
    \end{tabu}
\end{table}

\begin{figure}[ht]
    \center
    \caption{Merit of record gap for gaps up to 100,000}
    \includegraphics[scale=0.4]{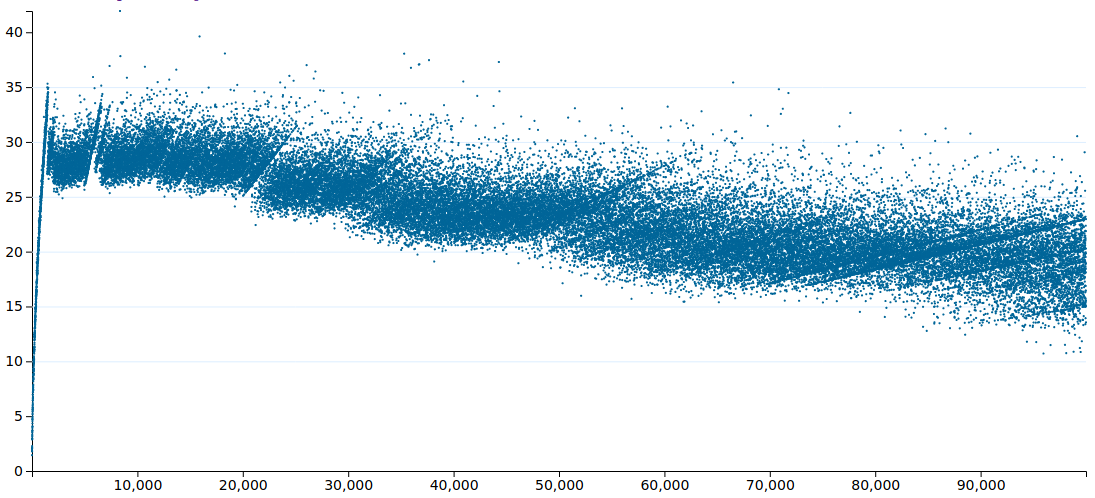}
    \label{figure:recordmerit}
\end{figure}

%\begin{figure}[ht]
%    \center
%    \includegraphics[scale=0.8]{5003_30030_100000000_80000_s38765_l50000M}
%    \caption{Stats and distribution for $M=100000000 + [0, 80000), X=38765, K=\frac{5003\#}{30030}$}
%    \label{appendix:probrecord}
%end{figure}

\begin{figure}[ht]
    \caption{Stats and distribution for $K=\frac{3001\#}{2190}\ \ M=1 + [0, 100000) $}
    \label{appendix:probrecord}
    \center
    \includegraphics[scale=0.5]{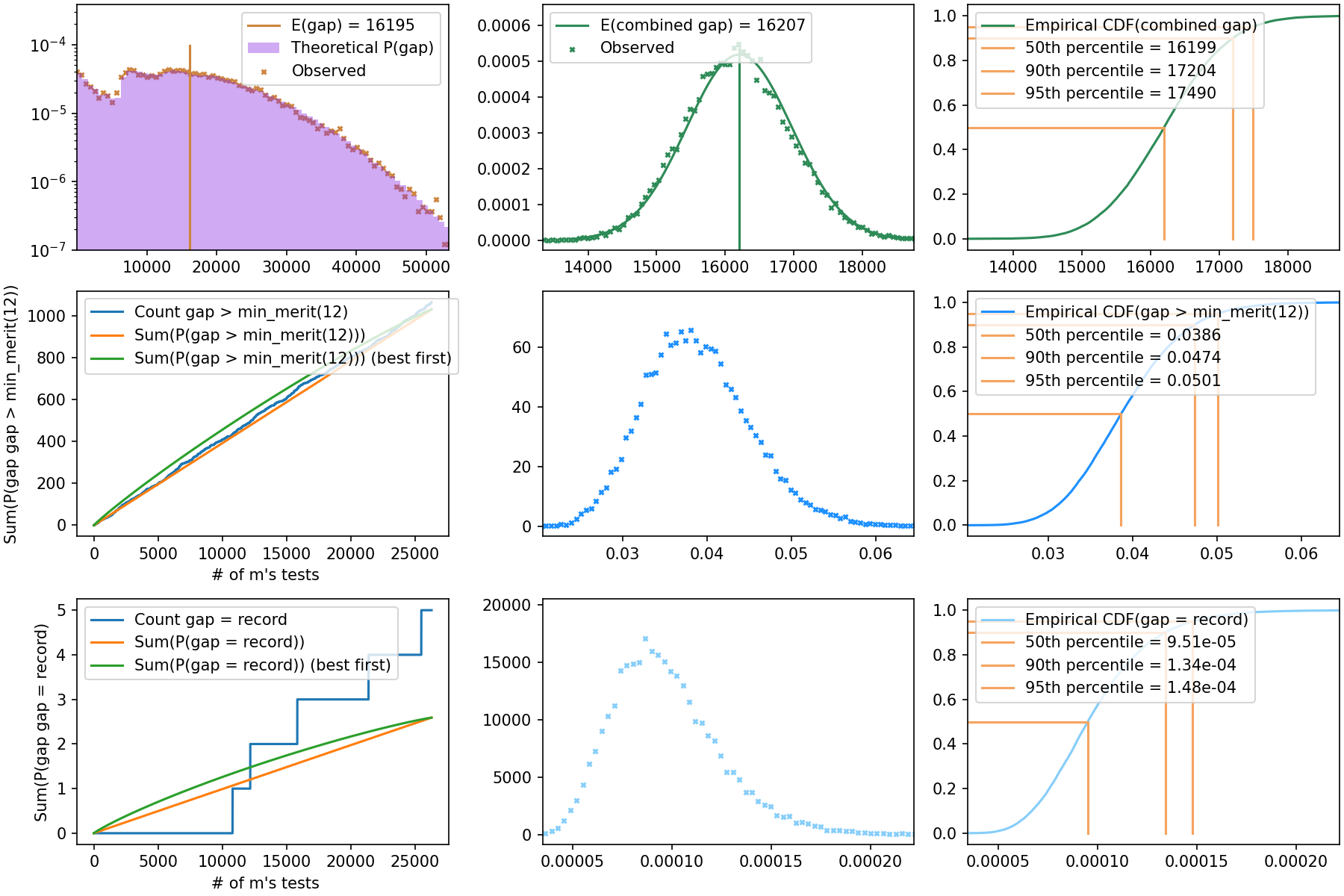} \\
    Top row: Probability of gap $X$, Distribution of expected gap, CDF of expected gap \\
    Middle row: Sum of probability of merit 12+ gap after $M$ tests, histogram of probabilities, CDF of probability \\
    Bottom row: Sum of probability of record gap after $M$ tests, histogram of probabilities, CDF of probability
\end{figure}

\end{document}